\documentclass[11pt]{article}

\usepackage[a4paper,margin=30mm]{geometry}
\usepackage[T1]{fontenc}
\usepackage{lmodern}
\usepackage{microtype}
\usepackage{mathtools,amssymb,amsthm}
\usepackage{booktabs,array,longtable}
\usepackage[shortlabels]{enumitem}
\usepackage{xcolor}
\usepackage[colorlinks=true,linkcolor=blue!55!black,citecolor=blue!55!black,urlcolor=blue!55!black]{hyperref}
\usepackage[nameinlink,noabbrev]{cleveref}
\usepackage{amsmath,amssymb,tikz}

\newtheorem{theorem}{Theorem}[section]
\newtheorem{proposition}[theorem]{Proposition}
\newtheorem{lemma}[theorem]{Lemma}
\theoremstyle{remark}

\newcommand{\R}{\mathbb R}
\newcommand{\Z}{\mathbb Z}
\newcommand{\one}{\mathbf 1}
\newcommand{\eps}{\varepsilon}
\newcommand{\Diag}{\operatorname{Diag}}

\title{A rigidity theorem for edge-connectivity under non-negative Bakry--\'Emery curvature}
\author{
 Kaizhe Chen\thanks{Department of Mathematics, Princeton University, Princeton, USA. Email: kc8556@princeton.edu}
 \and Shiping Liu\thanks{School of Mathematical Sciences, University of Science and Technology of China, Hefei, China. Email: spliu@ustc.edu.cn
  }\and Haotian Yang \thanks{Department of Mathematics, California Institute of Technology, Pasadena, USA. Email: hyang3@caltech.edu}
}
\date{}

\begin{document}
\maketitle
\begin{abstract}
The edge-connectivity \(\kappa(G)\) of a graph \(G\) is the minimum number of edges whose removal disconnects \(G\), and it is bounded above by the minimum degree \(\delta(G)\).
We classify all connected locally finite simple graphs with nonnegative Bakry--\'Emery curvature satisfying \(\kappa(G)<\delta(G)\).
All such graphs are infinite. Consequently, every finite connected simple graph with nonnegative Bakry--\'Emery curvature satisfies \(\kappa(G)=\delta(G)\),
confirming a conjecture of Chen, Koolen, and Liu. 
\end{abstract}

\section{Introduction}

Various notions of discrete Ricci curvature are defined in terms of the local structure of a graph \cite{NR17}. Recent work \cite{ChenKoolenLiu,ChenLiuYou,HPS24,LiuXia,ZZ26} has systematically explored the relationship between curvature bounds and global vertex- and edge-connectivity, providing insights into the structure of graphs with positive or nonnegative curvature.

The notion of Bakry--\'Emery curvature \cite{BE,Elworthy,LY10,Sch} is motivated by the Bochner identity and defined via \(\Gamma\)-calculus. On a graph, the curvature at a vertex \(x\) is determined entirely by the ball of radius \(2\) centered at \(x\) (see, e.g., \cite[Section 6]{CLP20}). Horn, Purcilly, and Stevens \cite{HPS24} established a lower bound for the vertex-connectivity of a connected graph in terms of its minimum degree and a positive lower bound on its Bakry--\'Emery  curvature. Subsequently, Chen, Koolen, and Liu \cite{ChenKoolenLiu} obtained the following lower bound for the edge-connectivity of connected graphs with nonnegative Bakry--\'Emery curvature.

\begin{theorem}[\cite{ChenKoolenLiu}]\label{thm:ckl}
If a connected locally finite simple graph has non-negative Bakry--\'Emery curvature, then \(\kappa(G)\geq\delta(G)-1\).
\end{theorem}

This lower bound on edge-connectivity is sharp. Chen, Koolen, and Liu \cite{ChenKoolenLiu} conjectured that every \emph{finite} connected simple graph with nonnegative Bakry--\'Emery curvature satisfies \(\kappa(G)=\delta(G)\). We confirm this conjecture by classifying all connected locally finite simple graphs with nonnegative Bakry--\'Emery curvature satisfying \(\kappa(G)=\delta(G)-1\). All such graphs are infinite.

\begin{theorem}\label{thm:main}
Let \(G\) be a connected locally-finite simple graph with non-negative Bakry--\'Emery curvature. Then \(\kappa(G)=\delta(G)\), unless \(G\) is isomorphic to one of
\[
 P_{\Z}\mathbin{\square}K_n\quad(n\geq 1),\qquad
 P_{\Z}^{\,2},\qquad P_{\Z}[K_2].
\]
Here \(P_{\Z}\) is the two-way infinite path; \(P_{\Z}\square K_n\) is the Cartesian product; \(P_{\Z}^{\,2}\) joins two integers precisely when their distance in \(P_{\Z}\) is one or two; and \(P_{\Z}[K_2]\) is the lexicographic product, so every \(P_{\Z}\)-layer is a \(K_2\) and consecutive layers span a \(K_{2,2}\). 
\end{theorem}



\begin{figure}[htbp]
\centering
\begin{tikzpicture}[
  x=1.55cm, y=1cm,
  vertex/.style={circle,fill=black,inner sep=0pt,minimum size=4.2pt},
  edge/.style={line width=0.6pt},
  fiber/.style={line width=0.95pt},
  every node/.style={font=\small}
]

\begin{scope}[yshift=3cm]

  \begin{scope}
    \clip (-2.60,-0.65) rectangle (2.60,0.65);

    \foreach \i in {-3,...,2} {
      \draw[edge] (\i, 0.50) -- ({\i+1}, 0.50);
      \draw[edge] (\i,-0.50) -- ({\i+1},-0.50);

      \draw[edge] (\i,-0.50) -- ({\i+1}, 0.50);
    }

    \foreach \i in {-2,...,2}
      \draw[fiber] (\i,-0.50) -- (\i,0.50);
  \end{scope}

  \foreach \i in {-2,...,2} {
    \pgfmathtruncatemacro{\evenlabel}{2*\i}
    \pgfmathtruncatemacro{\oddlabel}{2*\i+1}
    \node[vertex] at (\i, 0.50) {};
    \node[vertex] at (\i,-0.50) {};
    \node[above=5pt] at (\i, 0.50) {\(\evenlabel\)};
    \node[below=5pt] at (\i,-0.50) {\(\oddlabel\)};
  }
  \foreach \y in {-0.50,0.50} {
    \node at (-2.86,\y) {\(\cdots\)};
    \node at ( 2.86,\y) {\(\cdots\)};
  }
\end{scope}

\begin{scope}

  \begin{scope}
    \clip (-2.60,-0.65) rectangle (2.60,0.65);

    \foreach \i in {-3,...,2} {
      \draw[edge] (\i, 0.50) -- ({\i+1}, 0.50);
      \draw[edge] (\i,-0.50) -- ({\i+1},-0.50);
      \draw[edge] (\i, 0.50) -- ({\i+1},-0.50);
      \draw[edge] (\i,-0.50) -- ({\i+1}, 0.50);
    }

    \foreach \i in {-2,...,2}
      \draw[fiber] (\i,-0.50) -- (\i,0.50);
  \end{scope}

  \foreach \i in {-2,...,2} {
    \node[vertex] at (\i, 0.50) {};
    \node[vertex] at (\i,-0.50) {};
    \node[above=5pt] at (\i, 0.50) {\((\i,0)\)};
    \node[below=5pt] at (\i,-0.50) {\((\i,1)\)};
  }
  \foreach \y in {-0.50,0.50} {
    \node at (-2.86,\y) {\(\cdots\)};
    \node at ( 2.86,\y) {\(\cdots\)};
  }
\end{scope}

\end{tikzpicture}
\caption{The graphs \(P_{\mathbb Z}^{\,2}\) and
\(P_{\mathbb Z}[K_2]\). }
\label{fig:infinite-graphs}
\end{figure}
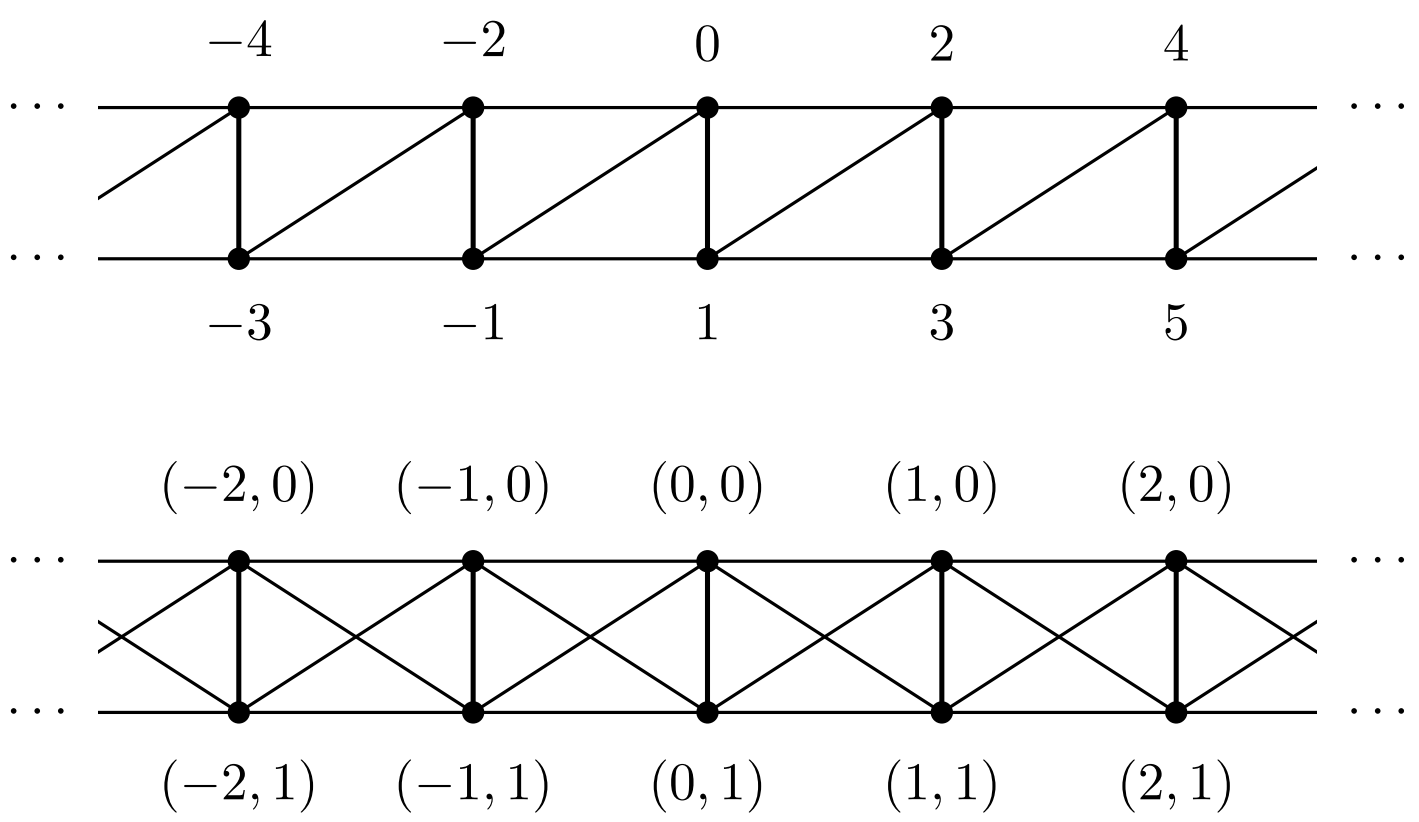

Our proof analyzes a minimum edge cut of a graph \(G\) with
nonnegative Bakry--\'Emery curvature and \(\kappa(G)<\delta(G)\).
Let \(t\) denote the maximum degree of the bipartite graph formed
by the cut edges. We show that \(t\in\{1,2\}\).

\textbf{Related work on Ollivier Ricci curvature.}
Ollivier Ricci curvature \cite{O09} and its Lin--Lu--Yau modification \cite{LLY} are
defined using optimal transport between probability measures
associated with vertex neighborhoods.
The study of these curvature notions and that of Bakry--\'Emery
curvature have often inspired one another.

Motivated in part by the results of Horn, Purcilly, and Stevens \cite{HPS24}
on Bakry--\'Emery curvature and vertex-connectivity, and by
Hehl's work \cite{Hehl25} on minimum degree and Lin--Lu--Yau curvature,
Chen, Liu, and You \cite{ChenLiuYou} systematically investigated the interaction between vertex-connectivity
and Lin--Lu--Yau curvature.
They also proved that every connected graph with positive Lin--Lu--Yau curvature satisfies \(\kappa(G)=\delta(G)\),
and asked whether the same conclusion holds for finite connected graphs with nonnegative Lin--Lu--Yau curvature. Liu and Xia \cite{LiuXia} answered this question affirmatively by classifying
all connected locally finite simple graphs with nonnegative Lin--Lu--Yau curvature satisfying \(\kappa(G)=\delta(G)-1\).

The present work is partly inspired by this classification.
All the exceptional graphs listed in \Cref{thm:main} also have
nonnegative Lin--Lu--Yau curvature.
Nevertheless, the classes of graphs with nonnegative
Bakry--\'Emery curvature and with nonnegative Lin--Lu--Yau
curvature are distinct \cite[Remark 2.10]{CKKLP2021}. 
Zhou and Zhang \cite{ZZ26} have also studied the relationship between
edge-connectivity and long-scale Ollivier and Lin--Lu--Yau curvatures.

\section{Preliminaries}

Throughout the paper, \(G=(V,E)\) is a connected locally-finite simple graph. For \(x\in V\), write \(N(x)\) for its neighborhood, \(d(x):=|N(x)|\), and \(\delta(G):=\min_{x\in V}d(x)\). For \(i\geq0\), let \(N_i(x)\) be the set of all vertices at distance $i$ to $x$. 
For $A\subseteq V$, put \(d_A(v)=|N(v)\cap A|\). Let $e(A)$ denote the number of edges with both endpoints in $A$.
For any disjoint \(A,B\subseteq V\), \(e(A,B)\) denotes the number of edges between $A$ and $B$. 
For a set \(S\subseteq V\), let \(\partial S :=E(S,V\setminus S)\) be the set of edges between $S$ and $V\setminus S$. 
A set \(S\) satisfying \(|\partial S|=\kappa(G)\) is called a \emph{shore} of a minimum cut.
Denote by \(L(G)\) the combinatorial Laplacian of \(G\).

For any integer $j>0$, we write \(P_j\) for a path with \(j\) edges and \(C_j\) for a cycle with \(j\) edges. For two graphs $G_1$ and $G_2$, the symbol \(G_1+G_2\) denotes their disjoint union. For any integer $m>0$, \(mP_1\) denotes the disjoint union of \(m\) copies of \(P_1\).

The (standard) graph Laplacian is defined as follows. For any function \( f:V\to\R\) and any vertex $x \in V$,
\[\Delta f(x):=\sum_{xy\in E}(f(y)-f(x)).\] 
For any two functions \(f,g:V\to\R\), define
\[
\Gamma(f,g):=\frac12\bigl(\Delta(fg)-f\Delta g-g\Delta f\bigr),\qquad
\Gamma_2(f,g):=\frac12\bigl(\Delta\Gamma(f,g)-\Gamma(f,\Delta g)-\Gamma(g,\Delta f)\bigr).
\]
We abbreviate \(\Gamma(f)=\Gamma(f,f)\) and \(\Gamma_2(f)=\Gamma_2(f,f)\). 
We say that the condition \(CD(\infty,0)\) holds at \(x\) if \(\Gamma_2(f)(x)\geq0\) for every \(f:V\to\R\). \(G\) has non-negative Bakry--\'Emery curvature if \(CD(\infty,0)\) holds at every vertex in $G$.

We write \(I\) and \(J\) for the identity and all-ones matrices of the appropriate order, respectively.  For any vertex set $A$, write \(\one_A\) for the indicator vector of a set \(A\).
We will use the curvature-matrix formulation of Bakry--\'Emery curvature developed in \cite{CKLP2022}, which characterizes the curvature-dimension condition in terms of positive-semidefiniteness of a matrix obtained by Schur complementation. See also the independent work of Siconolfi \cite{Si21}, and a related earlier work \cite{CLP20}. 

\begin{proposition}\label{prop:matrix}
Let $G=(V,E)$ be a graph. Let \(x\in V\) and \(d:=d(x)\). For each \(z\in N_2(x)\), put \(I_z:=N(z)\cap N(x)\) and \(k_z:=|I_z|\).
Then \(CD(\infty,0)\) holds at \(x\) if and only if
\(M_x\succeq0\), where
\begin{equation}\label{eq:Mx}
M_x:=\frac12J-\frac{d-3}{4}I+L(G[N(x)])+
\sum_{z\in N_2(x)}
\left(\frac34\Diag(\one_{I_z})-\frac1{k_z}\one_{I_z}\one_{I_z}^{\mathsf T}\right).
\end{equation}
Moreover, for any \(u=(u_y)_{yx\in E}\),
\begin{equation}\label{eq:quadratic}
u^{\mathsf T}M_xu=
\frac12\Bigl(\sum_{yx\in E}u_y\Bigr)^2+
\sum_{\{y,y'\}\subseteq N(x)}\omega_x(y,y')(u_y-u_{y'})^2
-\frac14\sum_{yx\in E}\eta_{xy}u_y^2,
\end{equation}
where
\[
\omega_x(y,y'):=\mathbf1_{yy'\in E}+
\sum_{\substack{z\in N_2(x)\\zy,zy'\in E}}\frac1{k_z}\ \ \ \ {\rm and}\ \ \ \
\eta_{xy}:=d(x)+d(y)-4-|N(x)\cap N(y)|.
\]
\end{proposition}

\begin{proof}
By \cite[Theorem~1.2]{CKLP2022}, \(CD(\infty,0)\) holds at \(x\) if and only if
\(A_\infty(x)\succeq 0\), where \(A_\infty(x)\) is the corresponding curvature matrix.
Moreover, with the normalization used here, \cite[Proposition~1.13(i)]{CKLP2022} gives \(A_\infty(x)=2M_x\).
So, \(CD(\infty,0)\) holds if and only if \(M_x\succeq0\).

To obtain \eqref{eq:quadratic}, observe that, for each
\(z\in N_2(x)\),
\[
\frac34\sum_{y\in I_z}u_y^2-\frac1{k_z}\left(\sum_{y\in I_z}u_y\right)^2=\frac1{k_z}\sum_{\{y,y'\}\subseteq I_z}(u_y-u_{y'})^2-\frac14 \sum_{y\in I_z}u_y^2.
\]
Together with
\[u^{\mathsf T}L(G[N(x)])u=\sum_{\{y,y'\}\in E(G[N(x)])}(u_y-u_{y'})^2,\]
this gives the coefficients \(\omega_x(y,y')\).
Moreover,
\[
\bigl|\{z\in N_2(x):y\in I_z\}\bigr|=d(y)-1-|N(x)\cap N(y)|.
\]
Thus, the remaining coefficient of \(u_y^2\) is \(-\eta_{xy}/4\), proving \eqref{eq:quadratic}.
\end{proof}

During the proof, we will use repeatedly the fact that a principal submatrix of a positive-semidefinite matrix is positive-semidefinite.

The following lemma is the case \(K=0\) of \cite[Lemma 3.1]{ChenKoolenLiu}, and the special case \(\varepsilon =0\) was first established in \cite{HPS24}.

\begin{lemma}\label{lem:cut-normal}
Let $G=(V,E)$ be a graph satisfying $CD(\infty,0)$ at $x\in V$. Suppose that $X,Y$ partition $N_1(x)$ and $A,B$ partition $N_2(x)$. Then, for any $\varepsilon\in \mathbb{R}$, we have
\begin{align*}
&(1-\eps)^2\left(e(X,Y)+e(Y,A)+e(X,B)
-\sum_{z\in A}\frac{d_Y(z)^2}{d_{N_1(x)}(z)}
-\sum_{z\in B}\frac{d_X(z)^2}{d_{N_1(x)}(z)}\right)\notag\\
\geq{}&\frac{d(x)-3}{4}\bigl(\eps^2|X|+|Y|\bigr)
+\frac14\bigl(\eps^2e(X,N_2(x))+e(Y,N_2(x))\bigr)
-\frac12\bigl(\eps|X|+|Y|\bigr)^2.
\end{align*}
\end{lemma}
The lemma follows by applying the curvature-dimension condition \(CD(\infty,0)\) at $x$, namely \(\Gamma_2(f)(x)\geq 0\), to suitably chosen test functions \(f\).

\section{Proof of the main result}

Assume, throughout this section, that a graph \(G=(V,E)\) satisfies \(CD(\infty,0)\) and \(\kappa(G)<\delta(G)\). For simplicity, write \(\delta:=\delta(G)\) and \(\kappa:=\kappa(G)\). 
By \Cref{thm:ckl}, \(\kappa=\delta-1\). 
Fix a minimum cut \(F=E(L,R)\) of size \(\kappa\). For any \(v\in L\), let \(h(v):=|N(v)\cap R|\), and define \(h\) symmetrically for $v\in R$. For any \(S\in\{L,R\}\) and finite proper subset \(S_1\subsetneq S\), put 
\[D_S(S_1):=\sum_{v\in S_1}(d(v)-2h(v))-2e(S_1).\]
Since $e(L,R)=\kappa$, we deduce that
\begin{align}\label{eq:defect}
    D_S(S_1)=e(S_1,S\setminus S_1)-e(S_1,V\setminus S)=|\partial(S\setminus S_1)|-\kappa\ge 0,
\end{align}
where equality holds if and only if \(S\setminus S_1\) is a shore.

\begin{lemma}\label{lem:defect}
Every shore is connected and has at least \(\delta(G)+1\) vertices.
\end{lemma}

\begin{proof}
If a shore $T$ has two components, their disjoint boundaries would each have at least \(\kappa\) edges. Altogether, the boundary of $T$ would have $2\kappa$ edges, a contradiction. So, every shore is connected.
Moreover, if a shore \(T\) had \(s\leq\delta\) vertices, then 
\[|\partial T|=\sum_{v\in T}d(v)-2e(T)\ge s\delta-s(s-1)=\delta+(s-1)(\delta-s)\ge\delta>\kappa,\]
which is a contradiction.
\end{proof}

Denote by \(V_F\) the set of vertices incident with at least one edge of \(F\).
Let \(H_F:=(V_F,F)\) be the bipartite graph consisting only of the cut edges, and let $t$ be the maximum degree of $H_F$. We next show that $t$ must be small.

Choose a vertex \(x\in V_F\) with \(h(x)=t\), and let \(S_x\in \{L,R\}\) be its shore. Put \(d:=d(x)\), $q:=d-t$. Notice that $q>0$ since, by  \Cref{lem:defect}, $S_x$ is connected. Set \(X:=N(x)\cap S_x\), and \(Y:=N(x)\setminus S_x\). Then, $|X|=q$ and $|Y|=t$. 
Take \(A:=N_2(x)\cap S_x\) and $B:=N_2(x)\setminus S_x$.
Set \[a:=e(X,Y),\quad b:=e(Y,A),\quad \text{and}\quad c:=e(X,B).\]
The local parameters \(d,q,a,b,c\) depend on the choice of \(x\), but we suppress this dependence in the notation. Since \(x\) is an arbitrary vertex with \(h(x)=t\), any constraints derived from this condition and the standing assumptions apply to the corresponding parameters at every such vertex, on either side of the cut.
Note that
\begin{equation}\label{mid1}
    \sum_{z\in A}\frac{d_Y(z)^2}{d_{N_1(x)}(z)}\ge \sum_{z\in A}\frac{d_Y(z)^2}{d_{Y}(z)+q}\ge \sum_{z\in A}\frac{d_Y(z)}{q+1}=\frac{b}{q+1}.
\end{equation}
Similarly, we deduce that
\begin{equation}\label{mid2}
    \sum_{z\in B}\frac{d_X(z)^2}{d_{N_1(x)}(z)}\ge \frac{c}{t+1}.
\end{equation}
Since every vertex in $X$ has at least $\delta -q$ neighbors in $Y\cup N_2(x)$, we find 
$$e(X,N_2(x))=e(X,Y\cup N_2(x))-a\ge q(\delta-q)-a.$$
So, \(e(X,N_2(x))\geq\max\{0,q(\delta-q)-a\}\), and similarly \( e(Y,N_2(x))\geq \max\{0, t(\delta-t)-a\}\). 
Substituting these two estimates and inequalities \eqref{mid1} and \eqref{mid2} into \Cref{lem:cut-normal} yields
\begin{align}\notag
&(1-\eps)^2\left(a+\frac{qb}{q+1}+\frac{tc}{t+1}\right)+\frac12(\eps q+t)^2\\ \label{eq:Px}
\geq {}&\frac{d-3}{4}(\eps^2 q+t)+\frac14\left(\eps^2\max\{0,q(\delta-q)-a\}+\max\{0,t(\delta-t)-a\}\right)
\end{align}
for every real \(\eps\). Since $e(L,R)=\kappa=\delta -1$, we have
\begin{equation}\label{abc}
a+\frac{qb}{q+1}+\frac{tc}{t+1}\le a+b+c\le \delta -t-1.
\end{equation}
By the definition of $t$, any vertex in $Y$ has at most $t$ neighbors in $S_x$ and hence at most $t-1$ neighbors in $X$. So, 
\begin{equation}\label{mid4}
    t(\delta-t)-a\ge t(\delta-t)-t(t-1)=t(\delta-2t+1).
\end{equation}
Thus, taking \(\eps=-t/q\) in \eqref{eq:Px} and using inequalities \eqref{abc} and \eqref{mid4} give
\begin{equation*}
\frac{(q+t)^2}{q^2}(\delta-1-t)\geq\frac{(d-3)t(q+t)}{4q}
+\frac{t }{4}\max\{0,\delta-2t+1\},
\end{equation*}
where we have discarded the non-negative term $\max\{0,q(\delta-q)-a\}$.
Since \(q=d-t\), we derive
\begin{equation}\label{eq:t-bound}
\delta-1-t\geq\frac{(d-3)t(d-t)}{4d}
+\frac{t(d-t)^2}{4d^2}\max\{0,\delta-2t+1\},
\end{equation}
Discarding the last term and using $d\ge \delta$ give
\begin{equation}\label{weak}
4d(d-1-t)\geq (d-3)t(d-t).
\end{equation}

Applying inequality \eqref{eq:defect} with \(S=S_x\) and \(S_1=\{x\}\) yields \( d\geq 2t\). 
It is direct to check that inequality \eqref{weak} does not hold for any positive integers $d,t$ with $d\ge 2t$ and $t\ge 5$. 
Hence, we have \(t\leq 4\). For \(t=4\), \eqref{eq:t-bound} leaves only \( (\delta,d)=(8,8)\) or \((9,9)\). 
If \( (\delta,d)=(8,8)\), then $q=4$ and inequality \eqref{abc} gives $a+b+c\le 3$. So, taking \(\eps=-1\) in \eqref{eq:Px} yields 
$$0\le \frac{9}{2}a+\frac{16}{5}(b+c)-18\le \frac{9}{2}(a+b+c)-18\le -\frac{9}{2},$$
which is a contradiction. Similarly, if \( (\delta,d)=(9,9)\), then taking \( \eps=-4/5\) in \eqref{eq:Px} leads to a contradiction. Therefore, we have \(t\leq 3\).

We shall use the following consequence of the Cauchy--Schwarz inequality. Let
\(U\subseteq N_2(x)\), and for \(z\in U\) put
\(v_z:=\one_{I_z}=\one_{N(z)\cap N(x)}\) and \(v:=\sum_{z\in U}v_z\).
Then
\[
\sum_{z\in U}\frac{v_zv_z^{\mathsf T}}{\one^{\mathsf T}v_z}\succeq
\frac{vv^{\mathsf T}}{\one^{\mathsf T}v}.
\]
Since \(\sum_{z\in U}\Diag(v_z)=\Diag(v)\), it follows that
\begin{equation}\label{cs}
    \sum_{z\in U}\left(\frac34\Diag(v_z)-\frac{v_zv_z^{\mathsf T}}{\one^{\mathsf T}v_z}\right)\preceq\frac34\Diag(v)-\frac{vv^{\mathsf T}}{\one^{\mathsf T}v}.
\end{equation}
Note that the left hand side appears in the expression for $M_x$ given in \eqref{eq:Mx}.

\subsection{Exclusion of \texorpdfstring{$t=3$}{t=3}}

In this section, our aim is to prove $t\le 2$. 
For the case \(t=3\), \eqref{eq:t-bound} leaves only \((\delta,d)=(6,6),(6,7),(7,7),\allowbreak (7,8),\allowbreak (8,8),\allowbreak (8,9),(9,9),(10,10)\), or \( (11,11)\). Write \eqref{eq:Px} as \(A_2\eps^2+A_1\eps+A_0\ge 0\) for any real $\varepsilon$, where
\[
A_2=\left( a+\frac{qb}{q+1}+\frac{3c}{4}\right)-\frac{(d-3)q+\max\{0,q(\delta-q)-a\}}4+\frac{q^2}{2}, 
\]
\[A_1=3q-2\left( a+\frac{qb}{q+1}+\frac{3c}{4}\right),\ \ \ A_0=\left( a+\frac{qb}{q+1}+\frac{3c}{4}\right)-\frac{3d+3\delta-a}{4}+9.\]
The conditions \(A_2\ge 0,A_0\geq 0\) and \(A_1^2\leq4A_2A_0\), together with inequality \eqref{abc}, leave only \((\delta,d,a,b,c)=(D,D,D-4,0,0)\), where \(6\leq D\leq10\).
Write
\[
X=\{x_1,\ldots,x_{D-3}\}
\qquad {\rm and}\qquad
Y=\{y_1,y_2,y_3\}.
\]

Since \(a=e(X,Y)=D-4\), there are exactly \(D-4\) edges between
\(X\) and \(Y\). For \(1\leq i\leq D-3\) and \(1\leq j\leq3\),
put
\[
\rho_i:=d_Y(x_i),\qquad
\sigma_j:=d_X(y_j),\qquad
r_i:=d_A(x_i),\qquad
s_j:=d_B(y_j).
\]
The equalities \(b=e(Y,A)=0\) and \(c=e(X,B)=0\) imply that every neighbor of \(x_i\) in \(N_2(x)\) lies in \(A\), and every neighbor of \(y_j\) in \(N_2(x)\) lies in \(B\).
Since \(d(x_i)\geq\delta= D\), while \(x_i\) has at most \(|X|-1=D-4\) neighbors in \(X\), we have \(r_i\geq3-\rho_i\). Similarly,
\(s_j\geq D-3-\sigma_j\). Hence, there are nonnegative integers \(\eta_i\) and \(\theta_j\) such that
\[r_i=3-\rho_i+\eta_i,\qquad s_j=D-3-\sigma_j+\theta_j.
\]

Put \(\alpha:=\sum_i\eta_i\) and
\(\beta:=\sum_j\theta_j\). Since
\(\sum_i\rho_i=\sum_j\sigma_j=e(X,Y)=D-4\), it follows that
\[e(X,A)=\sum_{1\leq i\leq D-3} r_i=2D-5+\alpha\qquad {\rm and} \qquad e(Y,B)=\sum_{1\le j\le 3} s_j=2D-5+\beta.\]
Applying \Cref{lem:cut-normal} with these exact edge counts gives
\[
\begin{split}
&(1-\eps)^2(D-4) \\
\geq{} &\frac{D-3}{4}\bigl(\eps^2(D-3)+3\bigr)
+\frac14\bigl((2D-5+\alpha)\eps^2+2D-5+\beta\bigr)
-\frac12\bigl(\eps(D-3)+3\bigr)^2
\end{split}
\]
for every \(\eps\in\R\). Considering this as a quadratic inequality in \(\eps\) leaves 
\[
\begin{array}{c|l}
D&(\alpha,\beta)\\ \hline
6&(0,0)\\
7&(0,0),(0,1),(1,0),(1,1),(2,0),(3,0)\\
8&(0,0),(0,1),(1,0),(1,1),(2,0),(2,1),
   (3,0),(4,0),(5,0)\\
9&(0,0),(0,1),(1,0),(2,0),(3,0),(4,0),(5,0),(6,0)\\
10&(0,0),(1,0),(2,0),(3,0),(4,0).
\end{array}
\]

Let \(H\) be the bipartite graph with vertex classes \(X\) and \(Y\)
whose edges are the edges of \(G\) between \(X\) and \(Y\).
Thus \(e(H)=D-4\), and \(\rho_i=d_H(x_i)\) and
\(\sigma_j=d_H(y_j)\) for any \(1\leq i\leq D-3\) and \(1\leq j\leq3\).
Let \(K_X\) and \(K_Y\) denote the complete graphs on \(X\) and
\(Y\), respectively, and write
\(r=(r_1,\ldots,r_{D-3})^{\mathsf T}\) and
\(s=(s_1,s_2,s_3)^{\mathsf T}\).
Adding to \(G[N(x)]\) any missing edges within \(X\) or within \(Y\)
adds positive-semidefinite Laplacians of edges to \(M_x\), see \cite[Proposition 6.10]{CLP20}.
Moreover, applying inequality \eqref{cs} separately to the vertices of \(A\) and \(B\) gives a further positive-semidefinite increase. Consequently,
\(M_x\preceq\widehat M\), where
\[
\widehat M=
\frac12J-\frac{D-3}{4}I
+L(K_X\cup K_Y\cup H)
+\frac34
\begin{pmatrix}
\Diag(r)&0\\
0&\Diag(s)
\end{pmatrix}
-
\begin{pmatrix}
\dfrac{rr^{\mathsf T}}{\sum_i r_i}&0\\[2mm]
0&\dfrac{ss^{\mathsf T}}{\sum_j s_j}
\end{pmatrix}.
\]

For every pair \((\alpha,\beta)\) in the preceding table, we consider every
bipartite graph \(H\) with \(e(H)=D-4\), up to relabeling the
vertices within \(X\) and within \(Y\), and every choice of
nonnegative integers $\eta_1+\cdots+\eta_{D-3}=\alpha$ and $\theta_1+\theta_2+\theta_3=\beta$.
This gives \(31\,344\) matrices \(\widehat M\). Computing their
determinants gives \(\det\widehat M<0\) in every case. Hence
\(\widehat M\not\succeq0\), which implies \( M\not\succeq0\), a contradiction.

\subsection{The case \texorpdfstring{$t=2$}{t=2}}

In this subsection, assume that $t=2$. Then $\delta=\kappa+1\ge 3$. 
If \(\delta=3\), the two
edges of \(F\) meet at the vertex \(x\). The inequality \eqref{eq:Px} requires
\[
R(\eps):=\frac{d-3}{4}(q\eps^2+2)+\frac{\eps^2\max\{0,q(3-q)\}+2}{4} -\frac{(q\eps+2)^2}{2}\leq 0
\]
for every \(\eps\in\mathbb R\). However, for \(d=3,4,5,6\), the choices \( \eps =-1,-2,-1,-1/10\) give \(R(\eps)=1/2,3,5/2,3/4\), respectively, while \(R(0)=(d-6)/2>0\) for \(d\geq7\). Thus, \(\delta\geq4\).
By $|Y|=t=2$ and the definition of $t$, we have
\begin{equation}\label{ab}
a+b=\sum_{y\in Y}(h(y)-1)\leq \sum_{y\in Y}(t-1)=2.
\end{equation}
Write \eqref{eq:Px} as \(A_2\eps^2+A_1\eps+A_0\ge 0\). Then $A_0\ge 0,A_2\ge 0$, and $A_1^2-4A_2A_0\le 0$. By \eqref{abc} and \eqref{ab},
\[
a+\frac{qb}{q+1}+\frac{2c}{3}\leq a+b+\frac{2c}{3}\leq\frac{2\delta-4}{3}.
\]
Therefore,
\[
A_0=a+\frac{qb}{q+1}+\frac{2c}{3}+2-\frac{d-3}{2}-\frac{2(\delta-2)-a}{4}\le \frac{11-\delta}{3}-\frac{d-\delta}{2},
\]
which implies \(\delta\leq 11\).
For \(4\leq\delta\leq 11\), combining $A_0\ge 0,A_2\ge 0$, \(A_1^2\leq4A_2A_0\), together with \eqref{abc} and \eqref{ab}, gives precisely
\begin{equation}\label{eq:t2-parameters}
\begin{array}{c|l}
\delta& (d,a,b,c)\\ \hline
4&(4,1,0,0),(5,1,0,0)\\
5&(5,0,1,1),(5,0,2,0),(5,1,0,1),(5,1,1,0),(5,2,0,0),
   (6,1,1,0),(6,2,0,0),(7,2,0,0)\\
6&(6,0,2,1),(6,1,0,2),(6,1,1,1),(6,2,0,0),(6,2,0,1),
   (7,1,1,1),(7,2,0,1)\\
7&(7,0,2,2),(7,1,0,3),(7,1,1,2),(7,2,0,1),(7,2,0,2),
   (8,2,0,2)\\
8&(8,1,1,3),(8,2,0,3)\\
9&(9,2,0,4).
\end{array}
\end{equation}

By the parameter table, we always have $a+b>0$, so $Y$ contains a vertex of degree two in \(H_F\).
By symmetry, every vertex of degree two in \(H_F\) has a neighbor of degree two in \(H_F\).
Therefore, every component of \(H_F\) is a single edge, a path with at least three edges, or an even cycle.

For \(S\in\{L,R\}\), put
\(n_2(S):=|\{w\in V_F\cap S:h(w)=2\}|\).
If \(H_F\) has a component consisting of a single edge,
say \(uv\) with \(u\in L\) and \(v\in R\), then
\(h(u)=h(v)=1\).
For \(w\in N(u)\setminus\{v\}\), we have
\(|N(w)\cap N(v)|\leq1+h(w)\).
Thus, \eqref{eq:Mx} gives
\begin{align}\notag
0\leq(M_v)_{uu}
&\leq\frac{6-d(u)-d(v)}4
+\sum_{w\in(V_F\cap L)\setminus\{u\}}
  \frac{h(w)}{1+h(w)}\\ \label{eq:t2-single-edge-diagonal}
&=\frac12-\frac{n_2(L)}3
-\frac{d(u)+d(v)-2\delta}{4}.
\end{align}
Consequently, \(n_2(L)\leq1\), and similarly \(n_2(R)\leq1\).
It follows that \(H_F\cong P_3+(\delta-4)P_1\).

For \(\delta=8,9\), the parameter table gives \(a+b=2\).
Since any vertex of degree two in \(H_F\) can be chosen as \(x\), the identity \(a+b=2\) holds with respect to each such choice, forcing both of its neighbors in \(H_F\) to have degree two.
Hence, every component of \(H_F\) other than a single edge
is a cycle. The preceding bound on \(n_2(L)\) and \(n_2(R)\)
then rules out single-edge components.
Using \(e(H_F)=\delta-1\), we obtain the following list:
\begin{equation}\label{eq:t2-supports}
\begin{array}{c|l}
\delta&H_F\text{ up to isomorphism}\\ \hline
4&P_3\\
5&P_4,\ P_3+P_1,\ C_4\\
6&P_5,\ P_3+2P_1\\
7&P_6,\ P_3+3P_1,\ P_3+P_3,\ C_6\\
9&C_8,\ C_4+C_4.
\end{array}
\end{equation}

The two cases with \(\delta=9\) are impossible.
Indeed, the parameter table gives
\((d,a,b,c)=(9,2,0,4)\) at every vertex of \(V_F\).
Since \(2+a+c=|F|\), any two distinct vertices of \(V_F\)
on the same side of the cut are adjacent.
There are two vertices in \(B\cap V_F\), each having two
neighbors in \(X\) and two in \(Y\), and
\(e(Y,B\setminus V_F)=8\).
For \(w:=\one_Y\in\mathbb R^{N(x)}\), \eqref{eq:Mx} gives
\[
w^{\mathsf T}M_xw
=2-3+2+2\left(\frac32-1\right)-\frac84=0.
\]
On the other hand, \(|X\setminus V_F|=4\), and
\((M_xw)_z=1\) for every \(z\in X\setminus V_F\),
since \(h(z)=0\) and \(b=0\).
This contradicts \(M_x\succeq0\), which would imply
\(M_xw=0\).

By \Cref{lem:t2-appendix}(i), the remaining cases in \eqref{eq:t2-supports} are $(\delta,H_F)\in\{(4,P_3),(5,P_4),(5,C_4)\}$.

\begin{lemma}\label{lem:no-P4}
The graph \(H_F\) is not isomorphic to \(P_4\).
\end{lemma}

\begin{proof}
Suppose that $H_F=p_0p_1p_2p_3p_4$, where $p_0,p_2,p_4\in L$ and $p_1,p_3\in R$.
Then \(\delta=5\). 
If \(p_1p_3\notin E\), then, in the notation used in \eqref{eq:Px} with \(x=p_1\), the two remaining edges \(p_2p_3,p_3p_4\) of \(F\) give
\((a,b,c)=(0,1,0)\), which does not occur in the \(\delta=5\) row of the parameter table.
Hence \(p_1p_3\in E\), and the corresponding row is \((d,a,b,c)=(5,1,0,1)\). Thus \(d(p_1)=5\), and the same argument gives \(d(p_3)=5\).

Let \(u,v\) be the two vertices in
\(N(p_1)\setminus\{p_0,p_2,p_3\}\), and define
\(w\in\mathbb R^{N(p_1)}\) by $w_{p_0}=w_{p_2}=-1$, $w_{p_3}=0$, and $w_u=w_v=1$.
Put \(k:=|N(p_4)\cap\{p_0,p_2\}|\).
Substituting $w$ into \eqref{eq:quadratic} and using
\eqref{cs} for the undetermined contributions from
\(N_2(p_1)\) gives $w^{\mathsf T}M_{p_1}w\leq -3/4$ when $k=1$ and $w^{\mathsf T}M_{p_1}w\leq -7/12$ when $k=2$.
Hence, \(k=0\). Repeating the same argument at \(p_3\) gives
\(|N(p_0)\cap\{p_2,p_4\}|=0\). Thus, \(p_0,p_2,p_4\) are pairwise nonadjacent.

Now, put \(m:=|(N(p_0)\cap L)\cap N(p_2)|\). The vertex \(p_1\) is isolated in
\(G[N(p_0)]\). Among the vertices of \(N_2(p_0)\) adjacent to \(p_1\), the vertex \(p_2\) is adjacent in \(N(p_0)\) precisely to \(\{p_1\}\cup(N(p_0)\cap L\cap N(p_2))\), whereas \(p_3,u,v\) are adjacent there only to \(p_1\).
So, \eqref{eq:Mx} gives
\[
\begin{aligned}
(M_{p_0})_{p_1p_1}
&=\frac12-\frac{d(p_0)-3}{4}
+\left(\frac34-\frac1{m+1}\right)
+3\left(\frac34-1\right)\\
&=\frac{5-d(p_0)}4-\frac1{m+1}<0,
\end{aligned}
\]
since \(d(p_0)\geq\delta=5\). This contradicts
\(M_{p_0}\succeq0\).
\end{proof}

\begin{lemma}\label{lem:delta-four}
If \(\delta=4\), then \(H_F\cong P_3\) and
\((d,a,b,c)=(4,1,0,0)\).
Moreover, \(d(v)=4\) for every \(v\in V_F\) with \(h(v)=2\),
and the two vertices of \(V_F\) in each of \(L\) and \(R\)
are adjacent.
\end{lemma}

\begin{proof}
The preceding tables give \(H_F\cong P_3\) and \((d,a,b,c)=(4,1,0,0)\) or \((5,1,0,0)\).
Suppose that \(d=5\), and label \(X=\{x_1,x_2,x_3\}\) and \(Y=\{y_1,y_2\}\) so that \(x_1y_1\) is the unique edge between \(X\) and \(Y\).
Put \(r_i:=d_A(x_i)\) and \(s_j:=d_B(y_j)\).
Since \(b=c=0\) and every vertex has degree at least four, there are nonnegative integers \(\eta_i,\theta_j\) such that $(r_1,r_2,r_3)=(\eta_1,1+\eta_2,1+\eta_3)$ and  $(s_1,s_2)=(1+\theta_1,2+\theta_2)$.
Write \(\alpha:=\eta_1+\eta_2+\eta_3\) and
\(\beta:=\theta_1+\theta_2\).
Then \(e(X,A)=2+\alpha\) and \(e(Y,B)=3+\beta\).
Substituting these edge counts into \Cref{lem:cut-normal} gives
\[
\frac{14-\alpha}{4}\eps^2+4\eps+\frac{5-\beta}{4}\geq 0
\qquad\text{for every }\eps\in\mathbb R.
\]
Consequently,
\(\alpha\leq14\), \(\beta\leq5\), and
\((14-\alpha)(5-\beta)\geq64\).
These conditions force \(\beta=0\) and
\(\alpha\in\{0,1\}\). Thus, $(\eta_1,\eta_2,\eta_3)
\in\{(0,0,0),(1,0,0),(0,1,0),(0,0,1)\}$ and $(s_1,s_2)=(1,2)$.

Let \(H\) be the graph on \(X\cup Y\) in which \(X\) and
\(Y\) each induce a complete graph and \(x_1y_1\) is the
only edge between them.
With coordinates ordered as
\((x_1,x_2,x_3,y_1,y_2)\), put
\(r=(r_1,r_2,r_3,0,0)^{\mathsf T}\) and
\(s=(0,0,0,s_1,s_2)^{\mathsf T}\).
Applying \eqref{cs} separately to \(A\) and \(B\), and
adding the missing edges within \(X\) and within \(Y\),
gives
\[
M_x\preceq\widehat M:=
\frac12J-\frac12I+L(H)
+\frac34\Diag(r+s)
-\frac{rr^{\mathsf T}}{r_1+r_2+r_3}
-\frac{ss^{\mathsf T}}{s_1+s_2}.
\]
For the four choices of \((\eta_1,\eta_2,\eta_3)\) listed
above, the matrices \(12\widehat M\) are respectively
\(Q_0,\allowbreak Q_1,\allowbreak Q_2,\allowbreak Q_3\), where
\[
Q_0=
\begin{pmatrix}
36&-6&-6&-6&6\\
-6&27&-12&6&6\\
-6&-12&27&6&6\\
-6&6&6&29&-14\\
6&6&6&-14&14
\end{pmatrix},
\qquad
Q_1=
\begin{pmatrix}
41&-10&-10&-6&6\\
-10&29&-10&6&6\\
-10&-10&29&6&6\\
-6&6&6&29&-14\\
6&6&6&-14&14
\end{pmatrix},
\]
\[
Q_2=
\begin{pmatrix}
36&-6&-6&-6&6\\
-6&26&-14&6&6\\
-6&-14&29&6&6\\
-6&6&6&29&-14\\
6&6&6&-14&14
\end{pmatrix},
\qquad
Q_3=
\begin{pmatrix}
36&-6&-6&-6&6\\
-6&29&-14&6&6\\
-6&-14&26&6&6\\
-6&6&6&29&-14\\
6&6&6&-14&14
\end{pmatrix}.
\]
However, $\det Q_0=-3761316$, $\det Q_1=-4270266$, and $\det Q_2=\det Q_3=-4455216$.
Thus, all four possibilities contradict \(M_x\succeq0\), and hence \(d=4\).
Since \(x\) was an arbitrary vertex with \(h(x)=2\), both vertices of degree two in \(H_F\) have degree four in \(G\).
Applying \(a=1\) at each of them gives the desired adjacencies.
\end{proof}

By \eqref{eq:t2-supports},
\Cref{lem:t2-appendix}(i), and \Cref{lem:no-P4},
the remaining possibilities are
\((\delta,H_F)=(4,P_3)\) and \((5,C_4)\).

Suppose first that \(H_F\cong C_4\).
Put \(F_L:=V_F\cap L\) and \(F_R=V_F\cap R\).
By \Cref{lem:t2-appendix}(ii), both sets induce \(K_2\), and their vertices have degree five.
The four edges from \(F_L\) to \(L\setminus F_L\) therefore form another minimum cut, say $F'$.
Each vertex of \(F_L\) has degree two in $H_{F'}$, so the preceding exclusions imply that $H_{F'}$ is again \(C_4\).
By \Cref{lem:t2-appendix}(ii), its two endpoints in \(L\setminus F_L\) form another adjacent pair of degree-five vertices. Iterating this argument on both shores gives a two-way infinite sequence of disjoint copies of \(K_2\), with all four edges between consecutive copies. That is, \( G\cong P_{\mathbb Z}[K_2]\).

Now, suppose that \(H_F\) is the path \(abcd\), with \(a,c\in L\) and \(b,d\in R\).
By \Cref{lem:delta-four},
\(ac,bd\in E\) and \(d(b)=d(c)=4\).
Let \(v\in R\setminus\{b,d\}\) be the remaining neighbor of \(b\).
Then
\(\partial(R\setminus\{b\})=\{cd,bd,bv\}\)
is a minimum cut whose cut graph is the path \(cdbv\).
Applying \Cref{lem:delta-four} to this cut gives
\(d(d)=4\) and \(dv\in E\).
The same argument on the other shore gives \(d(a)=4\).
Iteration on both shores, as in the \(C_4\) case, now yields
\(G\cong P_{\mathbb Z}^{\,2}\).

\subsection{The case \texorpdfstring{$t=1$}{t=1}}

Assume \(t=1\), and write \(F:=\{b_ic_i:1\leq i\leq \kappa\}\), where \(B:=\{b_1,\ldots,b_{\kappa}\}\subseteq L\) and \(C:=\{c_1,\ldots,c_{\kappa}\}\subseteq R\).
For any $1\leq i\leq \kappa$, set \(\eps_i:=d(b_i)-(\kappa+1)\) and \(\varphi_i:=d(c_i)-(\kappa+1)\). 
Since $\delta=\kappa+1$, all $\eps_i$ and $\varphi_i$ are nonnegative integers.
Put
$$E_B=\sum_{1\leq i\leq \kappa}\eps_i\ \ \  {\rm and}\ \ \  E_C=\sum_{1\leq i\leq \kappa}\varphi_i.$$

\begin{lemma}\label{lem:clique-ports}
Both sets \(B\) and \(C\) induce copies of \(K_{\kappa}\). Moreover, \( \eps_i\in\{0,1\}\), \(\varphi_i\in\{0,1\}\), and \(\eps_i+\varphi_i\leq1\) for any $1\leq i\leq \kappa$. 
\end{lemma}

\begin{proof}
Let $r_i$ be the number of indices $j\neq i$ such that both $b_ib_j$ and $c_ic_j$ are edges. As $F$ is a matching, applying \eqref{eq:Mx} at the root $c_i$ gives
\begin{equation}\label{Mbb}
    (M_{c_i})_{b_ib_i}=\frac{2r_i+6-d(b_i)-d(c_i)}{4}.
\end{equation}
Since $d(c_i)\geq\delta=\kappa+1$, while $c_i$ has only one neighbor in $L$ and at most $\kappa-1$ neighbors in $C$, it has a neighbor $w\in R\setminus C$. As $F$ is a matching, \eqref{eq:Mx} gives $(M_{c_i})_{b_iw}=1/2$. 
Since $M_{c_i}\succeq 0$, the diagonal entry $(M_{c_i})_{b_ib_i}$ cannot vanish; otherwise, the $b_iw$ principal $2\times2$ minor would have negative determinant. Hence, \eqref{Mbb} gives $d(b_i)+d(c_i)\leq 2r_i+5$.
By the fact that $r_i\leq\kappa-1$, we have
$$2\kappa+2=2\delta \leq d(b_i)+d(c_i)\leq 2r_i+5\leq 2\kappa+3.$$
This implies $r_i=\kappa-1$ for every $i$, and hence both $B$ and $C$ induce copies of $K_{\kappa}$. Moreover,
$d(b_i)+d(c_i)\leq2\kappa+3$ gives $\eps_i+\varphi_i \leq 1$.
In particular, we have $\eps_i,\varphi_i\in \{0,1\}$.
\end{proof}

Note that \Cref{lem:clique-ports} implies $E_B+E_C\le \kappa$.
Without loss of generality, assume that \(E_B\leq E_C\).  Then, 
\begin{equation}\label{Eb}
    E_B\leq\left\lfloor \frac{\kappa}{2} \right\rfloor.
\end{equation}
For any $1\leq i\leq \kappa$, put \(P_i:=N_{L\setminus B}(b_i)\), so \( |P_i|=1+\eps_i\).  
For any vertex \(z\in L\setminus B\), put \(X_z:=N_B(z)\), \(r_z:=|X_z|\), \( Y_z:=N_{L\setminus B}(z)\), and \(s_z:=|Y_z|\). Let 
\[x_z:=\sum_{b_i\in X_z}\eps_i ,\ \ \   q_z:=e(X_z,Y_z) ,\ \ \ {\rm and} \ \ \ \Omega_z:=\sum_{b\in X_z,\ y\in Y_z}\omega_z(b,y).\]
By \Cref{lem:clique-ports}, we deduce that 
\begin{equation}\label{xz}
    x_z=e(X_z,L\setminus (B\cup\{z\})).
\end{equation}
The pair \((r_z,x_z)\) is called the \emph{type} of \(z\). A vertex $z\in L\setminus B$ with \( r_z=\kappa\) is called a \emph{hub}.

\begin{lemma}\label{lem:port-conductance}
For any \(z\in L\setminus B\), \(s_z\geq1\), \(0\leq q_z\leq x_z\), and
\begin{equation}\label{eq:omega-bound}
\Omega_z\leq q_z+\frac{r_z(\kappa-r_z-E_B+x_z)}{r_z+1}
+\frac{2r_z(E_B-x_z)}{r_z+2}+\frac{(x_z-q_z)s_z}{s_z+1}.
\end{equation}
Moreover, we have $\Phi(r_z,s_z,x_z)\ge 0$ and $\Psi(r_z,s_z,x_z,q_z)\ge 0$, where
\begin{align}
\Phi(r_z,s_z,x_z)
&=\left(x_z+\frac{r_z(\kappa-r_z-{E_B}+x_z)}{r_z+1}+\frac{2r_z({E_B}-x_z)}{r_z+2}\right)\!\!\left(1+\frac{r_z}{s_z}\right)^2\notag\\
&\qquad-\frac14\left(r_z(\kappa+s_z-2)+
\frac{r_z^2}{s_z^2}(s_z(\kappa+r_z-2)-x_z)\right),\label{eq:Phi}\\
\Psi(r_z,s_z,x_z,q_z)
&=\left(q_z+\frac{r_z(\kappa-r_z-{E_B}+x_z)}{r_z+1}+\frac{2r_z({E_B}-x_z)}{r_z +2}+ \frac{(x_z-q_z)s_z}{s_z+1}\right)\!\!\left(1+\frac{r_z}{s_z}\right)^2\notag\\
&\qquad-\frac14\left(r_z(\kappa+s_z-2)+x_z-q_z+
\frac{r_z^2}{s_z^2}(s_z(\kappa+r_z-2)-q_z)\right).\label{eq:Psi}
\end{align}
\end{lemma}

\begin{proof}
Since $d(z)\ge\delta>\kappa$ and $r_z\le \kappa$, we have \(s_z\geq1\). Inequality \eqref{xz} yields \(q_z\leq x_z\). 
Next, we prove \eqref{eq:omega-bound}.
By the definition of $\omega_z$,
\begin{align*}
\Omega_z=q_z+\sum_{\substack{b\in X_z\\ y\in Y_z}}\
\sum_{\substack{w\in N_2(z)\\ w b,wy\in E}}\frac{1}{k_w}
&=q_z+\sum_{w\in N_2(z)}\frac{d_{X_z}(w)d_{Y_z}(w)}{k_w}\\
&=q_z+\sum_{w\in N_2(z)}\frac{d_{X_z}(w)d_{Y_z}(w)}{d_{X_z}(w)+d_{Y_z}(w)}.
\end{align*}
We now estimate the last sum according to the location of $w$. By \Cref{lem:clique-ports}, if $w=b_i\in B\setminus X_z$ and $\varepsilon_i=0$, then $d_{X_z}(w)=r_z$ and $d_{Y_z}(w)\leq1$, so its contribution is at most
$r_z/(r_z+1)$.
If $w=b_i\in B\setminus X_z$ and $\varepsilon_i=1$, then $d_{X_z}(w)=r_z$ and $d_{Y_z}(w)\leq2$, so its contribution is at most $2r_z/(r_z+2)$.
There are $\kappa-r_z-(E_B-x_z)$ and $E_B-x_z$ such vertices, respectively.
Finally, consider $w\in N_2(z)\setminus B$. If $d_{X_z}(w)$ or $d_{Y_z}(w)$ is 0, its contribution is 0. So, we may assume $w\in N_2(z)\setminus (B\cup R)$ and $d_{X_z}(w)\geq 1$. Since $d_{Y_z}(w)\leq s_z$, its contribution is at most $s_zd_{X_z}(w)/(s_z+1)$.
By \eqref{xz}, the sum of $d_{X_z}(w)$ over all such $w$ is at most $x_z-q_z$. 
Altogether, the inequality \eqref{eq:omega-bound} then follows.

In \eqref{eq:quadratic}, let $u$ be equal to $1$ on $X_z$ and $-r_z/s_z$ on $Y_z$. Then, the sum of $u_v$ over all neighbors $v$ of $z$ equals 0.
Since $u$ is constant on each of $X_z$ and $Y_z$, we calculate that
$$u^{\mathsf T}M_zu=\Omega_z\left(1+\frac{r_z}{s_z}\right)^2-\frac14\left(
\sum_{b\in X_z}\eta_{zb}+\frac{r_z^2}{s_z^2}\sum_{y\in Y_z}\eta_{zy}\right).$$
Since $B$ is a clique and $F$ is a matching, counting common neighbors gives
$$\sum_{b\in X_z}\eta_{zb}=r_z(\kappa+s_z-2)+x_z-q_z,\qquad\sum_{y\in Y_z}\eta_{zy}\geq s_z(\kappa+r_z-2)-q_z.$$
Together with $M_z\succeq0$, these estimates imply
$$0\leq\Omega_z\left(1+\frac{r_z}{s_z}\right)^2-\frac14\left(r_z(\kappa+s_z-2)+x_z-q_z+\frac{r_z^2}{s_z^2}\bigl(s_z(\kappa+r_z-2)-q_z\bigr)\right).$$
Substituting \eqref{eq:omega-bound} into the above inequality gives \eqref{eq:Psi}. Since $\Psi(r_z,s_z,x_z,q_z)$ is increasing in $q_z$ and $q_z\leq x_z$, replacing $q_z$ by $x_z$ gives the weaker inequality \eqref{eq:Phi}.
\end{proof}

\begin{lemma}\label{lem:double-row}
We have
\begin{equation}\label{inciden}
\sum_{z\in L\setminus B}r_z=\kappa+E_B\ \ \ \ {\rm and}\ \ \ \ 
\sum_{z\in L\setminus B}x_z=2E_B.
\end{equation}
Moreover, for any $1\leq i\leq \kappa$ with $\varepsilon_i=1$, we have
\begin{equation}\label{eq:double-row}
7\sum_{z\in P_i}r_z \geq 4\kappa+42-8E_B.
\end{equation}
\end{lemma}

\begin{proof}
Since $|P_j|=1+\varepsilon_j$ for any $1\leq j\leq \kappa$, double counting the edges between $B$ and $L\setminus B$ gives
$$
\sum_{z\in L\setminus B}r_z
=\sum_{j=1}^{\kappa}|P_j|
=\kappa+E_B.
$$
Moreover, since $\varepsilon_j^2=\varepsilon_j$ for any $1\leq j\leq \kappa$, we find
$$
\sum_{z\in L\setminus B}x_z
=\sum_{z\in L\setminus B}\sum_{b_j\in X_z}\varepsilon_j
=\sum_{j=1}^{\kappa}\varepsilon_j|P_j|
=\sum_{j=1}^{\kappa}\varepsilon_j(1+\varepsilon_j)
=2E_B.
$$

Now, fix $i$ with $\varepsilon_i=1$ and write $P_i=\{u,v\}$. Let $W_i$ be the total $\omega_{b_i}$-weight between $\{u,v\}$ and $B\setminus{b_i}$. 
The direct edges between these two sets have total weight $r_u+r_v-2$. For $a\in L\setminus(B\cup\{u,v\})$, let $t_a:=|N(a)\cap{u,v}|$. If $t_a>0$ and $a$ contributes to $W_i$, then the second-neighbor profile of $a$ at $b_i$ has size $r_a+t_a$. Exchanging the order of summation therefore gives
$$W_i=r_u+r_v-2+\sum_{t_a>0}\frac{t_ar_a}{r_a+t_a}.$$
Since $t_a\in\{1,2\}$ and $r_a\geq 1$, we have $t_ar_a/(r_a+t_a)\leq 2r_a/3$. It follows from \eqref{inciden} that
\begin{equation}\label{eq:Wi-upper}
W_i\leq r_u+r_v-2+\sum_{t_a>0}\frac{2r_a}{3}\leq r_u+r_v-2+\frac{2(\kappa+E_B-r_u-r_v)}{3}.
\end{equation}

We next obtain a lower bound for $W_i$. Since $\varepsilon_i=1$, \Cref{lem:clique-ports} gives $\varphi_i=0$, and hence $d(b_i)=\kappa+2$ and $d(c_i)=\kappa+1$. Write $\tau:=\mathbf 1_{uv\in E}$. Since $u,v\in L\setminus B$, they have no neighbors in $R$. We calculate that
$$
\eta_{b_i u}=\kappa+s_u-1-\tau,\qquad
\eta_{b_i v}=\kappa+s_v-1-\tau,\qquad {\rm and}\qquad
\eta_{b_i c_i}=2\kappa-1.
$$
Moreover, for every $j\neq i$, the unique second neighbor of $b_i$ adjacent to both $c_i$ and $b_j$ is $c_j$, so $\omega_{b_i}(c_i,b_j)=1/2$, whereas $\omega_{b_i}(c_i,u)=\omega_{b_i}(c_i,v)=0$.
Apply \eqref{eq:quadratic} at the root $b_i$ to the vector taking value $1$ on $u$ and $v$, value $0$ on $B\setminus{b_i}$, and value $-4$ on $c_i$. The preceding identities give
\begin{equation}\label{eq:Wi-lower}
W_i\geq\frac{2\kappa+6+s_u+s_v-2\tau}{4}.
\end{equation}
Since $\delta=\kappa+1$, $d(u)=r_u+s_u$ and $d(v)=r_v+s_v$, we have $s_u+s_v\geq 2\kappa+2-r_u-r_v$. Together with $\tau\leq1$, \eqref{eq:Wi-lower} implies $W_i\geq(4\kappa+6-r_u-r_v)/4$.
Combining this with \eqref{eq:Wi-upper} gives the desired inequality \eqref{eq:double-row}.
\end{proof}

\begin{lemma}\label{lem:large-matching}
We have $E_B=0$.
\end{lemma}

\begin{proof}
According to Appendix~\ref{app:matching}, it suffices to prove the case $\kappa\ge 13$.
We claim that every vertex $z\in L\setminus B$ satisfies $r_z\in\{0,1,2,\kappa\}$. Assume, to the contrary, that \(3\leq r_z\leq \kappa-1\). For brevity, write $r:=r_z$, $s:=s_z$, and $x:=x_z$.
By \Cref{lem:port-conductance}, we have $\Phi(r,s,x)\ge 0$. 
Regarding \(s\) as the variable, we may rewrite $\Phi(r,s,x)$ as 
\[\Phi(s)=A-\frac{r}{4}s+\frac{B}{s}+\frac{C}{s^2},\]
where $A$, $B$, and $C$ depend only on $\kappa$, $E_B$, $r$, and $x$.
By the assumption \(3\leq r\leq \kappa-1\),
\[
C=r^2\left(\frac{
E_Br^2+(3r+2)x+r(\kappa r+2\kappa -r^2-2r)}{(r+1)(r+2)}+\frac{x}{4} \right)>0.
\]
Set $s_0:=\kappa+1-r$. Then $r+s=d(z)\ge\delta=\kappa+1$ gives $s\ge s_0$. Note that 
\[
\frac{rs_0}{2}-A=\frac{2r^2(\kappa-2E_B)+(12r+8)(r-x)+r(\kappa(3r^2+3r-2)-2r^3-2r^2-8r-8)}{4(r+1)(r+2)}.
\]
By \eqref{Eb}, the fact \( x\leq r\), and the assumption \(3\leq r\leq \kappa-1\), we find the above expression positive, and hence \( A< rs_0/2\). If $s_*> s_0$ is an interior critical point of $\Phi$, then $\Phi'(s_*)=0$. That is,
\[
\frac{B}{s_*}=-\frac{rs_*}{4}-\frac{2C}{s_*^2}.
\]
Substituting this into $\Phi(s_*)$ gives
\[
\Phi(s_*)=A-\frac{rs_*}{2}-\frac{C}{s_*^2}<\frac{r(\kappa+1-r)}{2}-\frac{rs_*}{2} -\frac{C}{s_*^2}=\frac{r(s_0-s_*)}{2} -\frac{C}{s_*^2}\le 0.
\]
Therefore, $\Phi(r,s,x)\ge 0$ and $s\ge s_0$ imply that $\Phi(r,s_0,x)\geq 0$.

Set $N(E_B,x):=4(r+1)(r+2)(\kappa+1-r)^2\Phi(r,s_0,x)$. Then, $N(E_B,x)\ge 0$. Moreover, we calculate that
\[
N(E_B,x)=r(r+2)L(\kappa,r)+4r^2(\kappa+1)^2{E_B}+(4(\kappa+1)^2(3r+2)+r^2(r+1)(r+2))x,
\]
where
\[
L(\kappa,r)=2(1-r)\kappa^3+(4r^2-3r+5)\kappa^2-(4r^3-r^2+3r-4)\kappa+2r^4-2r^3-3r^2-2r+1.
\]
The coefficients of \({E_B}\) and \(x\) are positive, so inequality \eqref{Eb} and the fact $x\le \min\{r,E_B\}$ give $N(\lfloor \kappa/2\rfloor,\min\{r,\lfloor \kappa/2\rfloor\})\ge 0$.
However, a direct Mathematica computation shows that $N(\lfloor \kappa/2\rfloor,\min\{r,\lfloor \kappa/2\rfloor\})< 0$ for all integers \(\kappa\geq13\) and \(3\leq r\leq\kappa-1\).
It follows that every vertex $z\in L\setminus B$ satisfies $r_z\in\{0,1,2,\kappa\}$, as claimed.

Now, assume, for a contradiction, that $E_B>0$. Then there exists $1\leq i\leq \kappa$ with $\varepsilon_i=1$ (that is, $|P_i|=2$).
By \eqref{eq:double-row} and \eqref{Eb}, we have $\sum_{z\in P_i}r_z\ge 6$. Together with the above claim, this implies that there exists $z\in P_i$ with $r_z=\kappa$. 
By \eqref{inciden}, \(\sum_zr_z=\kappa+{E_B}<2\kappa\), so such a vertex $z$ is unique. Moreover, \(x_z={E_B}\). For simplicity, write \(s:=s_z\).  
Then, substituting $(r_z,s_z,x_z,q_z)=(\kappa,s,E_B,q_z)$ into \eqref{eq:Psi} and rearranging gives
\[
q_z-2s\geq\frac{sH(\kappa,s,E_B)}
{\kappa^2s+5\kappa^2+8\kappa s+s^3+5s^2},
\]
where $H(\kappa,s,E_B)$ equals
\[
(\kappa-2)s^3+(\kappa^2-\kappa-10)s^2
+\kappa(2\kappa^2-3\kappa-18)s
+2\kappa^2(\kappa-6)-\bigl(4\kappa^2+8\kappa s+3s^2-s\bigr)E_B.
\]
The coefficient of \(E_B\) in \(H\) is negative. Hence, by
\eqref{Eb},
\[
H(\kappa,s,E_B)\geq H\left(\kappa,s,\frac{\kappa}{2}\right)=(\kappa-2)s^3+(\kappa^2-\frac{5}{2}\kappa-10)s^2+(2\kappa^3-7\kappa^2-\frac{35}{2}\kappa)s-12\kappa^2.
\]
For all \(\kappa\geq 13\) and \(s\geq1\), the right-hand side is positive.
Therefore, \(q_z>2s\). On the other hand, since each \(y\in Y_z\) satisfies \( r_y\leq2\) by the claim, we derive
\[
q_z=e(B,Y_z)=\sum_{y\in Y_z}r_y\leq 2s,
\]
where we have used $X_z=B$ as $r_z=\kappa$. This leads to a contradiction.
\end{proof}

Since \(E_B=0\), every \(b_i\) has a unique neighbor \(\phi_B(i)\) in \(L\setminus B\). 

If a vertex $v\in L\setminus B$ satisfies \(|\phi_B^{-1}(v)|\ge 2\), set $r:=|\phi_B^{-1}(v)|$ and \( Y:=N_{L\setminus B}(v)\). 
Let $u$ be the vector equal to \(1\) on $\phi_B^{-1}(v)$ and \(-r/|Y|\) on \(Y\).
Since $N(v)=\phi_B^{-1}(v)\sqcup Y$, the chosen vector satisfies $\sum_{xv\in E}u_x=0$. Moreover, since $B$ is a clique, the total $\omega_v$-weight between $\phi_B^{-1}(v)$ and $Y$ is at most $r(\kappa-r)/(r+1)$. For $b_i\in\phi_B^{-1}(v)$, we have $\eta_{vb_i}=\kappa+|Y|-2$, while for $y\in Y$, $\eta_{vy}\ge \kappa+r-2.$
Therefore, \eqref{eq:quadratic} gives
\begin{equation}\label{uMu}
    \frac{u^{\mathsf T}M_vu}{r}\leq
    \frac{\kappa-r}{r+1}\left(1+\frac r{{|Y|}}\right)^2-
    \frac{\kappa+{|Y|}-2}{4}-\frac{r(\kappa+r-2)}{4{|Y|}}.
\end{equation}
Put $k:=\kappa-r$ and $a:={|Y|}-k-1$. Then $k\ge 0$, and since $d(v)\ge\delta=\kappa+1$, $a\geq 0$. Multiplying the right-hand side of \eqref{uMu} by $4(r+1){|Y|}^2$ yields the following cubic in $a$:
\begin{equation*}
\begin{aligned}
& -(r+1)a^3-\bigl(4kr+(r+1)^2\bigr)a^2 -\Bigl((5r-3)k^2+(3r-4)(r+1)k +(r+1)(2r^2-1)\Bigr)a \\
& -2(r-1)k^3-(r+1)(2r-5)k^2 -(r+1)(2r^2-3r-4)k-(r-1)(r+1)(2r+1).
\end{aligned}
\end{equation*}
When $r\geq 3$, the above expression is negative, and hence $u^{\mathsf T} M_vu<0$. When $r=2$, its coefficients $-3$, $-(8k+9)$, $-(7k^2+6k+21)$, and $-2k^3+3k^2+6k-15$ are all negative for any integer $k\ge 0$.
We also have $u^{\mathsf T}M_vu<0$, contradicting $M_v\succeq0$. Hence, $\phi_B$ is injective.

If \(E_C>0\), then there exists \(i\) with \(\varphi_{i}=1\). Put \( v_0:=\phi_B(i)\). Injectivity of $\phi_B$ gives \(d_{L\setminus B}(v)\geq \delta-1=\kappa\). Note that $(M_{b_i})_{c_ic_i}=1/4$, $(M_{b_i})_{c_iv_0}=1/2$, while \[(M_{b_i})_{v_0v_0}\leq 1+\frac{\kappa-1}{2}-\frac{\kappa+d_{L\setminus B}(v_0)}{4}\leq \frac{1}{2}.\]
Thus, the corresponding $2\times2$ principal minor is at most \(-1/8\), which is contradictory to $M_{b_i}\succeq0$. It follows that \(E_C=0\). By symmetry, there exists an injective map $\phi_C$ such that every \(c_i\) has a unique neighbor \(\phi_C(i)\) in \(R\setminus C\). 

Note that all the above arguments in this subsection hold for any minimum cut. Moreover, $\{c_i\phi_C(i):1\le i\le \kappa\}$ is a minimum cut with shores $L\cup C$ and $R\setminus C$. Therefore, $\phi_C(C)$ also induces a clique $K_{\kappa}$, and any vertex in $\phi_C(C)$ has a unique neighbor in $R\setminus (C\cup \phi_C(C))$ such that these neighbors are pairwise distinct. An iteration argument on both ends shows that $G \cong P_{\Z}\square K_n$.

\subsection{Verification of the exceptional graphs}
We have shown that $G$ is isomorphic to \(P_{\Z}\square K_n\), \(P_{\Z}^{\,2}\), or \(P_{\Z}[K_2]\). In this subsection, we prove that they indeed have nonnegative curvature.
All three families are vertex-transitive, so it suffices to check \eqref{eq:Mx} at one root. For \(P_{\Z}\square K_n\), order the neighbors as the \(n-1\) vertices in the root layer followed by the two axial neighbors. Then we have by \Cref{prop:matrix}
\[
M_x=\left(\frac{3n}{4}I_{n-1}-\frac12J_{n-1}\right)\oplus\frac12J_2.
\]
For \(n\geq2\), its eigenvalues are \(3n/4\) (with multiplicity \(n-2\)), \((n+2)/4\), \(1\), and \(0\); for \(n=1\), only the last block remains. Hence, \( M_x\succeq0\). We refer to \cite[Theorem 1.10]{CKLP2022} for a general result on the curvature matrices of Cartesian products of graphs. 

For \(P_{\Z}^{\,2}\) and \(P_{\Z}[K_2]\), we have by \Cref{prop:matrix}
\[
M_x=\frac14
\begin{pmatrix}
5&-4&2&2\\
-4&10&-2&2\\
2&-2&10&-4\\
2&2&-4&5
\end{pmatrix}
\qquad {\rm and}\qquad
M_x=
\begin{pmatrix}
4&-\frac12&-\frac12&-\frac12&-\frac12\\
-\frac12&\frac52&-\frac32&\frac12&\frac12\\
-\frac12&-\frac32&\frac52&\frac12&\frac12\\
-\frac12&\frac12&\frac12&\frac52&-\frac32\\
-\frac12&\frac12&\frac12&-\frac32&\frac52
\end{pmatrix},
\]
respectively. In both cases, \(M_x\succeq 0\), with smallest eigenvalue zero. This completes the proof of \Cref{thm:main}.

\section*{Acknowledgement}

S.L. is supported by the National Key R and D Program of China 2023YFA1010200 and the Scientific Research Innovation Capability Support Project for Young Faculty (Grant No. SRICSPYF-ZY2025160). Part of this paper was written during S.L.'s visit to the Simons Institute for the Theory of Computing at the University of California, Berkeley. S.L. gratefully acknowledges the Institute's hospitality and financial support.

\section*{Statement on AI}
The authors used AI as an exploratory tool to assist in developing the arguments. The authors
have carefully checked all mathematical proofs and take full responsibility for the content of the
manuscript.

\appendix

\section{Finite calculations for the \texorpdfstring{\(t=2\)}{t=2} case}\label{app:nonmatching}

\begin{lemma}\label{lem:t2-appendix}
Assume that \(t=2\). Then the following statements hold.
\begin{enumerate}[(i)]
\item \(H_F\) is not isomorphic to any of
\(P_3+P_1, P_5, P_3+2P_1, P_6, P_3+3P_1, P_3+P_3, C_6\).
\item If \(H_F\cong C_4\), then \(d(v)=5\) for every
\(v\in V_F\), and both \(G[V_F\cap L]\) and
\(G[V_F\cap R]\) are isomorphic to \(K_2\).
\end{enumerate}
\end{lemma}

\begin{proof}
We retain the notation used in \eqref{eq:Px}.
To prove~(i), suppose that \(H_F\) is one of the listed graphs.
The possible values of \((d,a,b,c)\) are listed in \eqref{eq:t2-parameters}.
We first describe the matrix bounds used to exclude these possibilities;
the remaining finite calculations are supplied in the supplementary verification files.

Taking \(u=\one\) in \eqref{eq:Mx} gives
$0\leq4\one^{\mathsf T}M_x\one =d^2+3d-e(N_1(x),N_2(x))$.
Since \(d\) is bounded by \eqref{eq:t2-parameters}, there are only finitely many choices, up to relabeling, for \(G[N(x)]\) and the incidences between \(N_1(x)\) and \( N_2(x)\).
The finite verification uses \eqref{eq:Mx} to construct the corresponding curvature matrices and \eqref{cs} to bound the contributions of vertices whose incidences are recorded only through their sums. Adding missing edges inside \(X\) gives a further positive-semidefinite increase.

For \(H_F\cong P_3+kP_1\), \(1\leq k\leq3\), let \(uv\) be a single-edge component, with \(u\in L\) and \(v\in R\).
Since \(n_2(L)=1\), \eqref{eq:t2-single-edge-diagonal} gives \(0\leq(M_v)_{uu}\leq
1/6-(d(u)+d(v)-2\delta)/4\).
Thus, \(d(u)=d(v)=\delta\).
Moreover, if a vertex of \((V_F\cap L)\setminus\{u\}\) were not adjacent to \(u\), its summand could be omitted, decreasing the upper bound by at least \(1/2\). Hence \(u\) is adjacent to every other vertex of \(V_F\cap L\), and the analogous assertion holds for \(v\).
Applying this argument to every single-edge component determines all same-shore adjacencies incident with its endpoints.
Only the two possible same-shore edges among the vertices of the \(P_3\) component remain to be chosen.

When both of these edges are present, \eqref{eq:t2-parameters} forces the two vertices with \(h=2\) to have degree \(\delta\). For \(k=1,2\), the supplementary verification files also apply \eqref{eq:Mx} and \eqref{cs} at an endpoint of a single-edge component of \(H_F\).

The supplementary verification files enumerate the remaining choices and give, for each of them, a matrix \(Q\) and a vector \(w\) such that \(w^{\mathsf T}Qw<0\).
Here \(Q\succeq M_v\) at the chosen root \(v\), or \(Q\succeq T^{\mathsf T}M_vT\), where the columns of \(T\) are the indicator vectors of the parts of a partition of \(N(v)\).
Thus \(T^{\mathsf T}M_vT\) represents the restriction to functions constant on each part. The files specify these partitions and prove that the enumerated choices cover every configuration permitted by the preceding bounds.
By \Cref{prop:matrix}, each negative quadratic value is a contradiction. This proves~(i).

For~(ii), write \(V_F\cap S_x=\{x,x'\}\) and \(Y=\{y_1,y_2\}\). Suppose that \(xx'\notin E\). Then \(x'\in N_2(x)\cap S_x\) and \(b=2\), so \eqref{eq:t2-parameters} gives \((d,a,b,c)=(5,0,2,0)\).
Put \(k=|N(x')\cap X|\) and \(\rho=\mathbf1_{y_1y_2\in E}\).
Thus \(0\leq k\leq3\) and \(\rho\in\{0,1\}\).
For \(\xi,\eta\in\mathbb R\), let \(w\) take the value \(\xi\) on \(Y\) and \(\eta\) on \(X\).
There are at least \(6-2\rho\) incidences from \(Y\) to \(N_2(x)\setminus V_F\), and at least \(6-k\) from \(X\) to that set. Hence \eqref{eq:Mx} gives
\[
 w^{\mathsf T}M_xw\leq
 \begin{pmatrix}\xi&\eta\end{pmatrix}
 C_{k,\rho}
 \begin{pmatrix}\xi\\\eta\end{pmatrix},
 \qquad
 C_{k,\rho}:=
 \begin{pmatrix}
 1+\rho/2-4/(k+2)&3-2k/(k+2)\\
 3-2k/(k+2)&3/2+2k/(k+2)
 \end{pmatrix}.
\]
However, \(\det C_{k,\rho}=((10+7\rho)k+6\rho-84)/(4(k+2))<0\) for all the stated values of \(k\) and \(\rho\).
This contradicts \Cref{prop:matrix}.
Thus \(xx'\in E\), and symmetry gives \(y_1y_2\in E\).
By \eqref{eq:t2-parameters}, \(d(v)\in\{5,6,7\}\) for every \(v\in V_F\).
Since \(x'\) is adjacent to \(x,y_1,y_2\), it has at most four other neighbors.
Using this bound together with \eqref{eq:Mx} and \eqref{cs}, the finite verification in the supplementary files excludes \(d(x)=6\) and \(d(x)=7\).
Consequently every vertex of \(V_F\) has degree five, proving~(ii).
\end{proof}

\section{The \texorpdfstring{$t=1$}{t=1} case for \texorpdfstring{\(\kappa\leq 12\)}{κ ≤ 12}}\label{app:matching}

In this section, we prove the case $\kappa\le 12$ of \Cref{lem:large-matching}.
Assume, for a contradiction, that \(E_B>0\). By \eqref{Eb}, we have
\(1\leq E_B\leq\lfloor\kappa/2\rfloor\).

We first show that \(s_z\le 199\) for any $z\in L\setminus B$ with $r_z>0$. Suppose, to the contrary, that \(s_z\geq 200\). Since
\(1\leq r_z\leq\kappa\), \(0\leq q_z\leq x_z\leq E_B\), and
\(E_B\leq 6\), the first factor in \eqref{eq:Psi} satisfies
\[
q_z+\frac{r_z(\kappa-r_z-E_B+x_z)}{r_z+1}
+\frac{2r_z(E_B-x_z)}{r_z+2}+\frac{(x_z-q_z)s_z}{s_z+1}
\leq \kappa-r_z+E_B\leq17.
\]
Hence its positive contribution is at most
\(17(1+12/200)^2<20\), whereas the first negative term is at least
\[
\frac{r_z(\kappa+s_z-2)}4
\geq\frac{\kappa+s_z-2}{4}
\geq50.
\]
Thus \(\Psi(r_z,s_z,x_z,q_z)<0\), a contradiction. So, every feasible value of
\(s_z\) is at most \(199\).

Let \(z_1,\ldots,z_m\) be the vertices of \(L\setminus B\) with
\(r_{z_i}>0\), and write \(r_i:=r_{z_i}\) and \(x_i:=x_{z_i}\).
By \eqref{inciden},
\[
\sum_{i=1}^m r_i=\kappa+E_B,
\qquad
\sum_{i=1}^m x_i=2E_B.
\]
For distinct \(i,j\), let
$\mu_{ij}:=\bigl|\{\ell:\varepsilon_\ell=1,\,P_\ell=\{z_i,z_j\}\}\bigr|$.
Then \(\mu_{ij}=\mu_{ji}\), \(\mu_{ii}=0\), and
\[
\sum_{j\neq i}\mu_{ij}=x_i.
\]
Moreover, if \(\mu_{ij}>0\), then \eqref{eq:double-row} gives $7(r_i+r_j)\geq4\kappa+42-8E_B$.
For every pair with \(\mu_{ij}>0\), put
\(a_{ij}:=\mathbf1_{z_i z_j\in E}\), and define
\[
q_i:=\sum_{j\neq i}\mu_{ij}a_{ij},
\qquad
c_i:=\bigl|\{j:\mu_{ij}>0,\ a_{ij}=1\}\bigr|.
\]
By construction, \(q_i=q_{z_i}\) and \(c_i\leq s_{z_i}\).
Since \(d(z_i)=r_i+s_{z_i}\geq\kappa+1\), every configuration
therefore gives integers \(s_i\) satisfying
\[
\max\{\kappa+1-r_i,c_i\}\leq s_i\leq199,
\qquad
\Psi(r_i,s_i,x_i,q_i)\geq0.
\]

These conditions leave a finite set of integer states. The verification
first enumerates the multisets \(\{(r_i,x_i)\}\), then the multiplicities
\((\mu_{ij})\) and the prescribed adjacencies \(a_{ij}\), discarding a
state whenever one of the corresponding sets of admissible values of
\(s_i\) is empty. This leaves \(2916\) states in total, occurring only
for
$(4,2),\allowbreak (5,2),\allowbreak (6,2),\allowbreak(6,3),\allowbreak (7,2),\allowbreak(7,3),\allowbreak (8,3),\allowbreak(8,4),(9,3),(9,4),\allowbreak (10,4),(10,5),\allowbreak (11,5),(12,6),$
where the pairs denote \((\kappa,E_B)\).

The remaining states are excluded by a second finite verification.
For each state, the program first tests curvature matrices obtained at
vertices \(z_i\in L\setminus B\) and at vertices \(b_\ell\in B\),
using compressed test spaces and Loewner upper matrices whenever the
unrecorded part of the radius-two neighborhood is irrelevant except
through aggregate incidence bounds. Only states not excluded by these
tests are refined by enumerating the remaining adjacencies among the
\(z_i\), the admissible values of \(s_i\), and the multiplicities of
vertices outside \(B\cup\{z_1,\ldots,z_m\}\) according to their
neighborhoods in \(\{z_1,\ldots,z_m\}\). No edges farther from the
chosen root are required.

For each terminal branch, the verification chooses as a root either
some \(b_\ell\in B\) or one of the vertices \(z_i\in L\setminus B\)
adjacent to \(B\). The neighbors of the chosen root are partitioned
into disjoint coordinate groups, with singleton groups for the
coordinates retained individually, and \(T\) denotes the matrix whose
columns are the indicator vectors of these groups. Thus
\(T^{\mathsf T}M_vT\) is the restriction of the curvature quadratic
form to vectors that are constant on each group. Using the remaining
incidence bounds, the program constructs a symmetric matrix \(Q\)
satisfying \(Q\succeq T^{\mathsf T}M_vT\), and finds a principal
submatrix of \(Q\) with negative determinant. Hence \(M_v\succeq0\)
is impossible on that branch. The verification exhausts every
refinement of every retained first-stage state and leaves no unresolved
branch.

The supplementary verification package contains the proof of the
matrix construction and exhaustiveness, the complete enumeration
program, all recorded matrices and coordinate labels, an independent
certificate checker, and the execution log.

\end{document}